\documentclass[11pt,psfig]{article}
\usepackage{amssymb,amsmath,amscd}
\newtheorem{lem}{Lemma}
\newtheorem{conj}{Conjecture}
\newtheorem{thm}{Theorem}
\newtheorem{cor}{Corollary}
\newtheorem{defn}{Definition}

\newcommand{\p}{\prime}

\newcommand{\commentout}[1]{}

\begin{document}
\title{Reconstruction and Higher Dimensional Geometry}
\author{Hongyu He \\
Department of Mathematics \\
Louisiana State University \\
email: hongyu@math.lsu.edu \\
}
\date{}
\maketitle
\abstract{Tutte proved that, if two graphs, both with more than two vertices, have the same
collection of vertex-deleted subgraphs, then the determinants of the two corresponding adjacency matrices
are the same. In this paper, we give a geometric proof of Tutte's theorem using vectors and angles.
We further study the lowest eigenspaces of these adjacency matrices.}
\section{Introduction}
Given the graph $G=\{V, E\}$, let $G_i$ be the graph obtained by deleting the $i-$th vertex $v_i$.
 Fix $n \geq 3$ from now on.
Let $G$ and $H$ be two graphs of $n$ vertices. The main conjecture in reconstruction theory,
states that if $G_i$ is isomorphic to $H_i$ for every $i$, then $G$ and $H$ are isomorphic (up to a reordering of $V$).
This conjecture is also known as the Ulam's conjecture.\\
\\
The reconstruction conjecture can be formulated in purely algebraic terms.
 Consider two $n \times n$ real symmetric matrices $A$ and $B$. Let $A_i$ and $B_i$ be the
 matrices obtaining by deleting the $i$-th row and $i$-th column of $A$ and $B$, respectively.
\begin{defn}
Let $\sigma_i$ be a $n-1$ by $n-1$ permutation matrix. Let $A$ and $B$ be two $n \times n$ real symmetric matrices.
We say that $A$ and $B$ are hypomorphic if there exists a set of $n-1 \times n-1$ permutation matrices
$$\{\sigma_1, \sigma_2, \ldots, \sigma_n \},$$
such that $B_i=\sigma_i A_i \sigma_i^t$ for every $i$. Put $\Sigma=\{\sigma_1, \sigma_2, \ldots,
\sigma_n \}$. We write
$B= \Sigma(A)$. $\Sigma$ is called a hypomorphism.
\end{defn}
The algebraic version of the reconstruction conjecture can be stated as follows.
\begin{conj}~\label{alg} Let $A$ and $B$ be two $n \times n$ symmetric matrices.
If there exists a hypomorphism $\Sigma$ such that $B=\Sigma(A)$, then there exists
a $n \times n$ permutation matrix $\tau$  such that $B=\tau A \tau^t$.
\end{conj}
We start by fixing some notations. If $M$ is a symmetric real matrix, then the eigenvalues of $M$ are real. We write
$$eigen(M)=(\lambda_1(M) \geq \lambda_2(M) \geq \ldots \geq \lambda_n(M)).$$
If $\alpha$ is an eigenvalue of $M$, we denote the corresponding eigenspace by $eigen_{\alpha}(M)$. 
Let $\bold 1_n$ be the $n$-dimensional row vector $(1,1,\ldots,
1)$. We may drop the subscript $n$ if it is implicit. Put $J =\bold
1^t \bold 1$. If $A$ and $B$ are hypomorphic, so are $A+t J$ and
$B+ t J$. 
\begin{thm}[Tutte]
Let $B$ and $A$ be two real $ n \times n$ symmetric matrices.
If $B$ and $A$ are hypomorphic then
$\det (B-\lambda I+ t J)=\det (A-\lambda I+ t J)$ for all $t, \lambda \in \mathbb R$. 
\end{thm}
In this paper, we will study the geometry related to Conjecture ~\ref{alg}. \commentout{ Suppose that $A$ is positive definite. Then $A$ represents the distance matrix of a $n-1$-simplex, $S(A)$. If $A$ and $B$ are hypomorphic, then the faces of $S(A)$ will be congruent to the faces of $S(B)$.
Conjecture 1 can then be translated as follows.
\begin{conj}
Let $S$ and $T$ be two $n-1$-simplices. Suppose that the faces of $S$ are congruent to the faces of $T$. Then $S$ and $T$ must be congruent.
\end{conj}
The requirement that $A$ and $B$ are positive definite is not essential because we can always choose a $\lambda$ such that $A+\lambda I$ and $B+\lambda I$ are both positive definite. It follows that Conjecture 2 is equivalent to Conjecture 1.\\
\\
By analyzing the geometry related to the reconstruction of $n-1$-simplices, we obtain the following theorem.}
Out main result can be stated as follows.
\begin{thm}[Main Theorem]
Let $B$ and $A$ be two real $ n \times n$ symmetric matrices.
Let $\Sigma$ be a hypomorphism such that $B=\Sigma(A)$. Let $t$ be a real number. Then there exists an open interval $T$ such that for $t \in T$ we have
\begin{enumerate}
\item $\lambda_n(A+ t J)=\lambda_n(B+t J)$;
\item
$eigen_{\lambda_n}(A+ t J)$ and $eigen_{\lambda_n}(B+t J)$ are both one dimensional;
\item
$eigen_{\lambda_n}(A+ t J)=eigen_{\lambda_n}(B+t J).$
\end{enumerate}
\end{thm}
A similar statement holds for the highest eigenspaces.\\
\\
Since the sets of majors of $A+t J$ and of $B+t J$ are the same, for every $t \in T$ and $ \lambda \in \mathbb R$,
\begin{equation}~\label{1}
\det(A+t J-\lambda I)-\det(B +t J-\lambda I)=\det(A+t J)-\det(B+t J).
\end{equation}
If $t \in T$, by taking $\lambda=\lambda_n(A+t J)$, we obtain
$$\det(A+t J)-\det(B+t J)=\det(A+t J-\lambda I)-\det(B+t J-\lambda I)=0.$$
Since the above statement is true for $t \in T $, $\det(A+t J)=\det(B+t J)$ for every $t$. By Equation. ~\ref{1}, we obtain $\det (B-\lambda I+ t J)=\det (A-\lambda I+ t J)$ for all $t, \lambda \in \mathbb R$. 
This is Tutte's theorem, which was proved using rank polynomials and Hamiltonian circuits. I should also mention that Kocay ~\cite{k} found a simpler way to deduce the reconstructibility of characteristic polynomials. \\
\\
Here is the content of this paper. We begin by presenting a
positive semidefinite matrix $A+\lambda I$ by $n$ vectors in
$\mathbb R^n$. We then interpret the reconstruction conjecture as
a generalization of a congruence theorem in Eulidean geometry.
Next we study the presentations of $A+\lambda I$ under the
perturbation by $t J$. We
define a norm of angles in higher dimensions and
establish a comparison theorem. Our comparison
theorem then forces hypomorphic matrices to have the same lowest
eigenvalue and eigenvector.\\
\\
I would like to thank the referee for his valuable comments.
\section{Notations}
Unless stated otherwise,
\begin{enumerate}
\item all linear spaces in this paper will be finite dimensional real Euclidean spaces;
\item all linear subspaces will be equipped with the induced Euclidean metric;
\item all vectors will be column vectors;
\item vectors are sometimes regarded as points in $\mathbb R^n$.
\end{enumerate}
Let $U=\{u_1, u_2, \ldots u_m\}$ be an ordered set of $m$ vectors in $\mathbb R^n$. $U$ is also interpreted as a $n \times m$ matrix.
\begin{enumerate}
\item Let $\mathrm{conv \ } U$ be the convex hull spanned by $U$, namely,
$$\{\sum_{i=1}^m \alpha_i u_i \mid \alpha_i \geq 0, \sum_{i=1}^m \alpha_i=1 \}.$$
\item Let $\mathrm{aff \ } U$ be the affine space spanned by $U$, namely,
$$\{\sum_{i=1}^m \alpha_i u_i \mid \sum_{i=1}^m \alpha_i=1 \}.$$
\item Let $\mathrm{span \ } U  $ be the linear span of $U$, namely,
$$\{\sum_{i=1}^m \alpha_i u_i \mid \alpha_i \in \mathbb R \}.$$
\end{enumerate}
Then $\mathrm{conv \ } U \subset \mathrm{aff \ } U \subset \mathrm{span \ } U $.\\
\\
Let $A$ be a matrix. We denote the $(i,j)$-th entry of $A$ by $a_{ij}$. We denote the transpose of $A$ by $A^t$.
Let ${{\mathbb R^+}^n}$ be the set of vectors with only positive coordinates.

\section{Geometric Interpretation}
Fix a standard Euclidean space $(\mathbb R^n, (,))$.
\begin{defn}
Let $A$ be a symmetric positive semidefinite real matrix. An ordered set of vectors $V=\{v_1,v_2, \ldots v_n \}$ is said to be a presentation of $A$ if and only if $(v_i, v_j)=a_{ij}$.
\end{defn}
Regarding $v_i$ as column vectors and $V$ as a $n \times n$ matrix, $V$ is a presentation of $A$ if and only if $V^t V= A$. Every positive semidefinite real matrix $A$ has a presentation. In addition, the presentation $V$ is
unique up to a left multiplication by an orthogonal matrix.
\begin{defn}
Let $S$ and $T$ be two sets of
vectors in $\mathbb R^n$. $S$ and $T$ are said to be congruent
if there exists an orthogonal linear transformation
 in $\mathbb R^n$ that maps $S$
onto $T$.
\end{defn}
So $A=\sigma B \sigma^t$ for some permutation $\sigma$ if and only if $A$ and $B$ are presented by two congruent subsets in $\mathbb R^n$.\\
\\
Now consider two hypomorphic matrices $B=\Sigma(A)$. Observe that $B+\lambda I= \Sigma(A+\lambda I)$. Without loss of generality, assume $A$
and $B$ are both positive semidefinite. Let $U$ and $V$ be their presentations respectively.
Since $B_i=\sigma_i A_i \sigma_i^t$, $U-\{u_i\}$ is congruent to $V-\{v_i\}$. Then the reconstruction conjecture can be stated as follows.
\begin{conj}[Geometric reconstruction]
Let
$$S=\{u_1,u_2, \ldots, u_n \}$$
and
$$T=\{v_1,v_2, \ldots v_n \}$$
 be two finite
sets of vectors in $\mathbb R^m$. Assume that $S-\{u_i\}$ is congruent
 to $T-\{v_i\}$ for every $i$. Then $S$ and $T$ are congruent.
\end{conj}
Generically, $m=n$. 
\begin{defn}
We say that $U=\{u_i\}_1^n$ is in good
 position if the point $0$ is in the interior of the convex hull of $U$ and the convex hull of $U$ is of dimension $n-1$.
\end{defn}
\begin{lem}~\label{tfae} Let $A$ be a symmetric positive semidefinite matrix. The following are equivalent.
\begin{enumerate}
\item $A$ has a presentation in good position.
\item Every presentation of $A$ is in good position.
\item $rank(A)=n-1$ and $eigen_{0}(A) = \mathbb R \alpha$ for some $\alpha \in (\mathbb R^+)^n$.
\end{enumerate}
\end{lem}
Proof: Since $A$ is symmetric positive semidefinite, $A$ has a presentation. Let $U$ be a presentation of $A$. \\
\\
If $U$ is in good position, then every presentation obtained from an orthogonal linear transformation is also in good position.
Since a presentation is unique up to an orthogonal linear transformation,   $(1) \leftrightarrow (2)$. \\
\\
Suppose $U$ is in good position. Then $rank(U)=n-1$. So $rank(A)=n-1$. Since $0$ is in the interior of the convex hull of $U$, there exists $\alpha=(\alpha_1, \alpha_2, \ldots \alpha_n)^t$ such that
$$0=\sum_{1}^{n} \alpha_i u_i; \qquad \sum_{1}^n \alpha_i=1; \qquad \alpha_i >0 \,\, \forall \ i.$$
Since $rank(U)=n-1$, $\alpha$ is unique. Now $U \alpha= 0$ implies
$$A \alpha= U^t U \alpha= U^t 0= 0.$$
Since $rank(A)=rank(U)=n-1$, $eigen_0(A)=\mathbb R \alpha$. So $(2) \rightarrow (3)$. \\
\\
Conversely, suppose $rank(A)=n-1$ and $eigen_{0}(A)=\mathbb R \alpha$ with $\alpha \in {\mathbb R^+}^n$.
 Then $\sum_{i} \alpha_i u_i =0$ and the linear span $\mathrm{span \ } U$ is of dimension $n-1$. Thus, $0$ is in  $\mathrm{conv \ } U$. It follows that
$\mathrm{aff \ } U=\mathrm{span \ } U$. So $\dim(\mathrm{conv \ } U)=\dim(\mathrm{aff \ } U)=\dim(\mathrm{span \ } U)=n-1$. So $(3) \rightarrow (1)$. Q.E.D.

\begin{lem}~\label{alpha} Let $U$ be a presentation of $A$. Suppose that $U$ is in good position. Let $\alpha_i$ be the volume of the convex hull of $\{0, u_1, u_2, \ldots, \hat{u_i}, \ldots, u_n \}$. Then $$\alpha=(\alpha_1, \alpha_2, \ldots, \alpha_n)^t$$
is a lowest eigenvector.
\end{lem}
The proof can be found in many places. For the sake of completeness, I will give a proof using the language of exterior product. \\
\\
Proof: Choosing an orthonormal basis properly, we may assume that every $u_i \in \mathbb R^{n-1}$. $U$ becomes a $(n-1) \times n$ matrix. Let $x_1,x_2 \ldots x_{n-1}$ be the row vectors of $U$.  Consider the exterior product
$$x_1 \wedge x_2 \wedge \ldots \wedge x_{n-1}.$$
Let $\beta_i$ be the $i$-th coordinate in terms of the standard basis
$$\{(-1)^{i-1} e_1 \wedge e_2 \wedge \ldots \wedge \hat{e_i} \wedge \ldots \wedge e_n \mid i \in [1, n] \}.$$
Put $\beta=(\beta_1,  \beta_2, \ldots, \beta_n)^t$.
Notice that $x_i \wedge (x_1 \wedge x_2 \wedge \ldots \wedge x_{n-1})=0$ for $1 \leq i \leq n-1$. Therefore,
 $(x_i,\beta)=0$ for every $i$. So $U \beta=0$. It follows that $\sum_{i=1}^n  \beta_i u_i =0$.
Since $0$ is in the convex hull of $\{u_i\}_1^n$, $ \beta_i$ must be either all negative or all positive. Clearly,
$$|\beta_i|=|u_1 \wedge \ldots \wedge \hat{u_i} \wedge \ldots \wedge u_n|=(n-1)! \alpha_i.$$
Therefore, we have $U \alpha=0$. Then $A \alpha= U^t U \alpha=0$. $\alpha$ is a lowest eigenvector. Q.E.D.
\begin{thm}~\label{main2}
Suppose that $B=\Sigma(A)$. Suppose that $A$ and $B$ have presentations in good position. Then $eigen_0(A)=eigen_0(B) \cong \mathbb R$.
\end{thm}
Proof: Let $U$ and $V$ be presentations of $A$ and $B$ respectively. Then $U$ and $V$ are in good position.
Notice that the volume of the convex hull of
$$\{0, u_1,u_2, \ldots, \hat{u_i}, \ldots u_n \}$$
 equals the volume of the convex hull of
$$\{0, v_1,v_2, \ldots, \hat{v_i}, \ldots v_n \}$$
By Lemma ~\ref{alpha} and Lemma ~\ref{tfae}, $eigen_0( A)=eigen_0(B) \cong \mathbb R$. So the lowest eigenspace of $A$ is equal to the lowest eigenspace of $B$. Q.E.D.

\section{Perturbation by $J$}
Recall that $J=\bold 1_n^t \bold 1_n$.
We know that $B=\Sigma(A)$ if and only if $B+t J =\Sigma(A+t J)$. Let us see how presentations of $A+t J$ depend on $t$.
Let $A$ be a positive definite matrix.
Let $U=\{u_i\}_1^n$ be a presentation of $A$.\\
\\
Let
$\mathrm{aff \ } U$ be the affine space spanned by $U$.
Then $\{ u_i \}$ are affinely independent.
Let $u_0$ be the orthogonal projection of the origin onto $\mathrm{aff \ } U$. Then $(u_0, u_i-u_0)=0$ for every $i$. We obtain
$$ U^t u_0=\|u_0\|^2 \bold 1.$$
It follows that $u_0= \|u_0\|^2 (U^t)^{-1} \bold 1$. Consequently,
$$\|u_0\|^2=(u_0, u_0)=\|u_0\|^4 \bold 1^t U^{-1} (U^t)^{-1} \bold 1=\|u_0\|^4 \bold 1^t A^{-1} \bold 1.$$
Clearly, $\|u_0\|^2=\frac{1}{\bold 1^t A^{-1} \bold 1}$. We obtain the following lemma.
\begin{lem}~\label{u0}
Let $A$ be a positive definite matrix.
Let $U=\{u_i\}_1^n$ be a presentation of $A$. Let $u_0$ be the orthogonal projection of the origin onto $\mathrm{aff \ } U$. Then $\|u_0\|^2= \frac{1}{\bold 1^t A^{-1} \bold 1}$ and
$$u_0=\frac{1}{\bold 1^t A^{-1} \bold 1} (U^t)^{-1} \bold 1.$$
\end{lem}
Consider $\{u_i-s u_0\}_1^n$. Notice that
$$(u_i-s u_0, u_j-s u_0)=(u_i-u_0+(1-s)u_0, u_j-u_0+(1-s)u_0)=(u_i-u_0, u_j-u_0)+(1-s)^2 (u_0, u_0).$$
Taking $s=0$, we have
$$(u_i, u_j)=(u_i-u_0, u_j-u_0)+(u_0, u_0).$$
Therefore
$$(u_i-s u_0, u_j-s u_0)=(u_i,u_j)-(u_0,u_0)+(1-s)^2(u_0,u_0)=(u_i,u_j)+(s^2-2s)\|u_0\|^2.$$
We see clearly that $A+(s^2-2s) \|u_0\|^2 J$ is presented by $\{u_i-su_0\}_1^n$. Observe that $$\mathrm{span}(u_1- s u_0, u_2- s u_0, \ldots, u_n-s u_0)$$
is of dimension $n$ for all $s \neq 1$. So $A+(s^2-2s) \|u_0\|^2 J$ is positive definite for all $s \neq 1$. If $s=1$, we see that
$A-\|u_0\|^2 J$ is presented by $\{u_i-u_0\}_1^n$ whose linear span is of dimension $n-1$. We obtain the following lemma.
\begin{lem}~\label{per} Let $A$ be a symmetric positive definite matrix. Let $U$ be a presentation of $A$. Let $u_0$ be the orthogonal projection of the origin onto $\mathrm{aff \ } U$. Then $\{u_i - s u_0 \}_1^n$ is a presentation of $A+ (s^2-2s) \| u_0 \|^2 J$. Let $t=(s^2- 2s) \| u_0 \|^2$. Then $A + t J$ is positive definite for all $ t > - \| u_0 \|^2$ and positive semidefinite for $t=- \| u_0 \|^2$.
\end{lem}
Notice that 
$$u_0=\frac{1}{\bold 1^t A^{-1} \bold 1} (U^t)^{-1} \bold 1=\frac{1}{\bold 1^t A^{-1} \bold 1}U (U^{-1} (U^t)^{-1}) \bold 1=\frac{1}{\bold 1^t A^{-1} \bold 1}U A^{-1} \bold 1.$$
\begin{thm}~\label{per1}
Let $A$ be a symmetric positive definite matrix. Let $U$ be a presentation of $A$. Let $u_0$ be the orthogonal projection of the origin onto $\mathrm{aff \ } U$. Then $u_0=\frac{1}{\bold 1^t A^{-1} \bold 1}U A^{-1} \bold 1$ and  the following are equivalent.
\begin{enumerate}
\item $A -\|u_0 \|^2 J$ has a presentation in good position;
\item  $u_0$ is in the interior of $\mathrm{conv} \ U$;
\item $A^{-1} \bold 1 \in { \mathbb R^+ }^n$.
\end{enumerate}
\end{thm}
\begin{cor}~\label{lambda} Let $A$ be a real symmetric matrix.
There exists $\lambda_0 $ such that for every $\lambda \geq \lambda_0$ there exists a real number $t$ such that $A+\lambda I + t J$ has a presentation in good position.
\end{cor}
Proof: 
\commentout{
Write $A= P^t D P$ where $P$ is orthogonal and $D$ is diagonal. First choose $\lambda_1 \geq \| A \|$. Suppose that $\lambda \geq \lambda_1$. Then $D+ \lambda I$ is a diagnoal matrix with positive entries.  Let $(D+\lambda I)^{\frac{1}{2}}$ be the square root of $D+\lambda I$ with only positive entries. Let $U(\lambda)=  (D+\lambda I)^{\frac{1}{2}} P$. Then  $U(\lambda)$ is a presentation of $A+\lambda I$. Because of Theorem ~\ref{per1}, it suffices to show that there exists a $\lambda_0 \geq \lambda_1$ such that for every $\lambda \geq \lambda_0$, the orthogonal projection of the origin onto $\mathrm{aff \ } U(\lambda)$ is in the interior of $\mathrm{conv \ }U(\lambda)$.}
Instead, consider $I+s A$ with $s=\frac{1}{\lambda}$. $I+ s A$ is related to $A +\lambda I$ by a constant multiplication:
$$\lambda(I+ s A)=\lambda I+A. $$
Let $s_0=\frac{1}{\| A \| +1}$ where $\|A\|$ denote the operator norm.  Suppose that $0 \leq s \leq s_0$. Then $I+ s A$ is positive definite. For $s=0$, $(I+s A)^{-1} \bold 1 \in {\mathbb R^+}^n$. Since
$$s \rightarrow (I+s A)^{-1} \bold 1$$
is continuous on $(0, s_0)$, there exists a $s_1 \in (0, s_0)$ such that
$(I+s A)^{-1} \bold 1 \in {\mathbb R^+}^n$ for every $s \in (0,s_1]$.
So for every $\lambda \in [\frac{1}{s_1}, \infty)$, $(A+ \lambda I)^{-1} \bold 1=\lambda^{-1}(s A+ I)^{-1} \bold 1 \in  {\mathbb R^+}^n$.
Let $\lambda_0=\frac{1}{s_1}$. So for every $\lambda \geq \lambda_0$, $(A+ \lambda I)^{-1} \bold 1 \in {\mathbb R^+}^n$.
By Theorem ~\ref{per1}, for every $\lambda  \geq \lambda_0  $ there exists a $t$ such that 
$A+\lambda I + t J$ has a presentation in good position.
 Q.E.D. 
\section{Higher Dimensional Angle and Comparison Theorem}
\begin{defn}
Let $U=\{u_1,u_2,\ldots u_n\}$ be a subset in $\mathbb R^n$. $\mathbb R^n$ may be contained in some other Euclidean space. Let $u$ be a point in $\mathbb R^n$. The angle
$\angle(u, U)$ is defined to be the region
$$\{ \sum_{1}^{n} \alpha_i (u_i-u) \mid \alpha_i \geq 0 \}.$$
Two angles are congruent if there exists an isometry that maps one angle to the other.
Let $\mathcal B$ be the unit ball in $\mathbb R^n$. The norm of $\angle(u,U)$ is defined to be the volume of $\angle(u, U) \cap \mathcal B$, denote it by $|\angle(u, U)|$.
\end{defn}
Let me make a few remarks.
\begin{enumerate}
\item Firstly,
if two angles are congruent, their norms are the same. But, unlike the 2 dimensional case, if the norms of two angles are the same, these two angles may not be congruent.
\item Secondly, if $\{u_i -u \}_1^n$ are linearly dependent, then $|\angle(u, U)|=0$. If $u$ happens to be in $\mathrm{aff \ } U$, then $|\angle(u, U)|=0$.
\item
According to our definition, $|\angle(u, U)|$ is always less than half of the volume of $\mathcal B$.
\item More generally, one can allow $\{\alpha_i\}_1^n$ to be in a collection of other sign patterns which correspond to quadrants in two dimensional case. Then the norm of an angle can be greater than half of the volume of $\mathcal B$.
\end{enumerate}
\begin{lem}
If $\angle(u, U) \subseteq \angle(u, V)$, then $|\angle(u, U)| \leq |\angle(u, V)|$. If $|\angle(u, U)| >0$ and $\angle(u, U)
$ is a proper subset of $\angle(u, V)$ then $|\angle(u, U)| < |\angle(u, V)|$.
\end{lem}
\begin{thm}[Comparison Theorem]
Let $\angle(u, U)$ be an angle and $|\angle(u, U)| \neq 0$. Suppose that $v$ is contained
in the interior of the convex hull of $\{ u \} \cup U$. Then $|\angle(u, U)| < |\angle(v, U)|$.
\end{thm}
Proof: Without loss of generality, assume $u=0$. Suppose $|\angle(0, U)| >0$. Let $U=\{u_1, u_2, \ldots u_n \}$. Then $U$ is linearly independent. Since $v$ is in the interior of $ \mathrm{conv}(0,U)$, $v$ can be written as
$$ \sum_{i=1}^n \alpha_i u_i $$
with $\alpha_i \in \mathbb R^+$ and $\sum_{i}^n \alpha_i < 1$. \\
\\
Let $U^{\prime}=\{u_i-v \}_1^n$. It suffices to prove that $\angle(0, U)$ is a proper subset of  $\angle(0, U^{\p})$.
Let $x$ be a point in $ \angle(0, U)$ with $x \neq 0$. Then $x=\sum_{i} x_i u_i$ for some $ x_i \geq 0$ with $\sum_i x_i > 0$. Define for each $i$
$$y_i=x_i+ \alpha_i \frac{\sum x_i}{1-\sum \alpha_i}$$
The reader can easily verify that $\sum y_i(u_i-v)=x$.
\commentout{
\begin{equation}
\begin{split}
 \sum y_i(u_i-v) = & \sum [x_i+ \alpha_i \frac{\sum x_i}{1-\sum \alpha_i}](u_i-v) \\
=& \sum [x_i+ \alpha_i \frac{\sum x_i}{1-\sum \alpha_i}] u_i- (\sum [x_i+ \alpha_i \frac{\sum x_i}{1-\sum \alpha_i}])v \\
=& \sum x_i u_i+ [\frac{\sum x_i}{1-\sum \alpha_i}] [\sum \alpha_i u_i]-[\sum x_i] v- [\sum \alpha_i][ \frac{\sum x_i}{1-\sum \alpha_i}] v \\
=& x+[\frac{\sum x_i}{1-\sum \alpha_i}] v-[\sum x_i] v- [\sum \alpha_i][ \frac{\sum x_i}{1-\sum \alpha_i}] v \\
= & x
\end{split}
\end{equation}}
Observe that $y_i > x_i \geq  0$. So $\angle(0, U)$ is a proper subset of $\angle(0, U^{\prime})$. It follows that
$$|\angle(0, U)| < | \angle(0, U^{\prime})| = |\angle(v, U)|.$$
Q.E.D.

\begin{thm}~\label{com}
Let $U=\{u_1,u_2, \ldots, u_n \} \subset \mathbb R^m$ for some $m \geq n$. Suppose that $|\angle(u, U)| >0$. Suppose that the orthogonal projection of $u$ onto $\mathrm{aff \ } U$ is in the interior of $\mathrm{conv \ } U$. Let $v$ be a vector
such that $(u-v, u-u_i)=0$ for every $i$. If $u \neq v$ then  $|\angle(u,U)| > |\angle(v, U)|$.
\end{thm}
Proof: Without loss of generality, assume that $u=0$. Then $U$ is linearly independent and 
$v \perp u_i$ for every $i$. Let $u_0$ be the orthogonal projection of $u$ onto $\mathrm{aff \ } U$. By our assumption,
$u_0$ is in the interior of $\mathrm{conv \ } U$ and $\|u_0\| \neq 0$. Let 
$$v^{\prime}=(1-\sqrt{\frac{\|v\|^2}{\|u_0\|^2}+1}) u_0.$$
Then 
$$\|v^{\prime}-u_0 \|^2= (\frac{\|v\|^2}{\|u_0\|^2}+1) \|u_0\|^2= \|v\|^2+ \|u_0\|^2.$$
Notice that $v \perp u_i$ and   $u_0 \perp u_i-u_0$. We obtain
\begin{equation}
\begin{split}
(u_i-v^{\prime}, u_j-v^{\prime})= & (u_i-u_0+ \sqrt{\frac{\|v\|^2}{\|u_0\|^2}+1} u_0, u_j-u_0+ \sqrt{\frac{\|v\|^2}{\|u_0\|^2}+1} u_0) \\
= & (u_i-u_0, u_j -u_0)+(\frac{\|v\|^2}{\|u_0\|^2}+1) \|u_0\|^2 \\
= & (u_i-u_0, u_j-u_0)+(u_0, u_0)+(v, v)\\
= & (u_i, u_j)+(v,v)\\
= & (u_i-v, u_j-v).
\end{split}
\end{equation}
Hence $\angle(v, U) \cong \angle(v^{\prime}, U)$. Notice that 
$1-\sqrt{\frac{\|v\|^2}{\|u_0\|^2}+1} < 0$. So the origin sits between $v^{\prime}$ and 
$u_0$ which is in the interior of $\mathrm{conv \ } U$. Therefore, $0$ is in the interior of 
$\angle(v^{\prime}, U)$. By the Comparison Theorem,
$| \angle(v^{\prime}, U) | > | \angle(0, U)|$. Consequently,
$| \angle(v, U) | > | \angle(0, U)|$. Q.E.D.

\begin{thm}~\label{main1}
Suppose that $B=\Sigma(A)$. Let $\lambda_0$ be as in Cor. ~\ref{lambda} for both $B$ and $A$. Fix $\lambda \geq \lambda_0$. Let $t_1$ and $t_2$ be two real numbers such that $A+ \lambda I + t_1 J$ and $B+ \lambda I+t_2 J$ have presentations in good position. Then $t_1=t_2$.
\end{thm}
Proof: We prove by contradiction. Without loss of generality, suppose that $ t_1 > t_2$.
Let $U$  be a presentation of $A+\lambda I+t_1 J$. Then $U$ is in  good position. So $0$ is in the interior of $\mathrm{conv \ } U$.
 Let $V$ be a representation of $B+\lambda I + t_2 J$. Then $V$ is in good position. So $0$ is in the interior of $\mathrm{conv \ } V$ and $\dim (\mathrm{span \ } V)=n-1$. Let $v_0 \perp \mathrm{span \ } V $ and $\|v_0\|^2=t_1-t_2$. Let $V^{\prime}=\{v_i+v_0 \}_1^n$. Clearly, $V^{\prime}$ is a presentation of $B+\lambda I+ t_1 J$.\\
\\
By Thm. ~\ref{com}, for every $i$,
$$|\angle(0, V \backslash \{v_i\})| > |\angle(-v_0, V \backslash \{v_i\})| =|\angle(0, V^{\prime} \backslash \{v_i+v_0 \})|.$$
Since $B+\lambda I+ t_1 J=\Sigma(A+\lambda I+t_1 J )$, $V^{\prime} \backslash \{v_i+v_0\}$ is congruent to $U \backslash \{u_i\}$ for every $i$. Therefore $|\angle(0, V^{\prime} \backslash \{v_i+v_0 \})|=|\angle(0, U \backslash \{u_i\})|$. Since $0$ is in the interiors of the convex hulls of $U$ and of $V$, we have
$$Vol(\mathcal B)=\sum_{i=1}^n |\angle(0, V \backslash \{v_i\})|
> \sum_{i=1}^n |\angle(0, V^{\prime} \backslash \{v_i+v_0\})|=\sum_{i=1}^n |\angle(0, U \backslash \{u_i\})| =Vol(\mathcal B).$$
This is a contradiction. Therefore, $t_1=t_2$. Q.E.D.
\section{Proof of the Main Theorem}
Suppose $B=\Sigma(A)$. Suppose $\lambda_0$ satisfies Cor. ~\ref{lambda} for both $A$ and $B$. So for every $\lambda \geq \lambda_0$ there exist real numbers $t_1$ and $t_2$ such that $A+\lambda I + t_1 J$ has a presentation in good position and $B+\lambda I+ t_2 J$ has a presentation in good position. By Theorem ~\ref{main1}, $t_1=t_2$. Because of the dependence on $\lambda$, put $t(\lambda)=t_1=t_2$. By Theorem ~\ref{main2}, 
$$eigen_0(A+\lambda I+ t(\lambda) J)=eigen_0(B+\lambda I+ t(\lambda) J) \cong \mathbb R.$$
Since $0$ is the lowest eigenvalue of $A+ \lambda I+ t(\lambda) J$ and $B+\lambda I+ t(\lambda) J$, $\lambda$ is the lowest eigenvalue of $A +t(\lambda) J$ and $B+ t(\lambda) J$. In addition, 
$$eigen_{-\lambda}(A + t(\lambda) J)=eigen_{-\lambda}(B+t(\lambda) J) \cong \mathbb R.$$ 
Now it suffices to show that $t([\lambda_0, \infty))$ covers a nonempty open interval.\\
\\
By Lemme ~\ref{per} and Lemma ~\ref{u0},
$$t(\lambda)= -\| u_0\|^2= -\frac{1}{\bold 1^t (A+ \lambda I)^{-1} \bold 1}.$$
So $t(\lambda)$ is a rational function.
Clearly,  $t([\lambda_0, \infty))$ contains a nonempty open interval $T$. For $t \in T$, we have
$\lambda_n(A+ t J)= \lambda_n(B+ t J)$ and
$eigen_{\lambda_n}(A+ t J)=eigen_{\lambda_n}(B+ t J) \cong \mathbb R$.
This finishes the proof of Theorem 1. Q.E.D. \\
\\
\commentout{
 Since $\det(A+ t J)$ depends on $t$ continuous, $eigen(A+t J)$ depends on $t$ continuously. Therefore, there exists an interval $J$ containing $0$ such that $\lambda_{n-1}(A+t J)> \lambda_{n}(A+t J)$ for every $t \in J$.
So for $t \in J$, $eigen_{\lambda_n}(A+t J)$ and $eigen_{\lambda_n}(B+t J)$ are one dimensional. Notice that for $t \in J$, $$
t \rightarrow eigen_{\lambda_n}(A+t J) \in \mathbb RP^{n-1}$$ and
$$ t \rightarrow eiegn_{\lambda_n}(B+t J)$$
are both continuous. Therefore, there exists an interval $I$ such that
for $t \in I$, the lowest eigenvectors of $B+t J$ and $A+t J$ are in ${\mathbb R^{+}}^n$. By Lemma 1, for $t \in I$,
$B+t J-\lambda_n(B+t J) I$ and $A+t J-\lambda_n(A+t J)I $ both have presentations in good position. By Theorem ~\ref{main2},
$eigen_{\lambda_n}(A+t J)=eigen_{\lambda_n}(B+t J)$.}
Tutte's proof involves certain polynomials associated with a graph. It is algebraic in nature. The main instrument in our proof is the comparison theorem. Presumably, there is a connection between the geometry in this paper and the polynomials defined in
Tutte's paper. In particular, given $n$ unit vectors $u_1, u_2, \ldots u_n$, can we compute  the function $|\angle(0, U)|$ explicitly in terms of $U^t U=A$? This question turns out to be hard to answer. The norm $| \angle(0,U) |$ as a function of $A$ may be closely related to the functions studied in Tutte's paper ~\cite{tu}.

\end{document}